\documentclass[10pt,reqno]{amsart}

\usepackage{amssymb}
\usepackage{xypic}
\xyoption{all}
\def\umono{\ar@{_{(}->}[u]}
\def\uumono{\ar@{_{(}->}[uu]}

\def\lmono{\ar@{_{(}->}[l]}
\def\llmono{\ar@{_{(}->}[ll]}

\newcommand{\Z}{{\mathbb Z}}

\newcommand{\F}{{\mathbb F}}

\newcommand{\aut}[2]{\operatorname{aut}_{#1}(#2)}
\newcommand{\Aut}[2]{\operatorname{Aut}_{#1}(#2)}

\newcommand{\map}[3]{\operatorname{map}(#2,#3)_{#1}}
\newcommand{\Map}[5]{\operatorname{map}(#2,#3;#4,#5)_{#1}}
\newcommand{\func}[3]{\mbox{$#1 \colon #2 \to #3$}}

\newcommand{\A}{\ifmmode{\mathcal{A}}\else${\mathcal{A}}$\fi}
\newcommand{\K}{\ifmmode{\mathcal{K}}\else${\mathcal{K}}$\fi}
\newcommand{\U}{\ifmmode{\mathcal{U}}\else${\mathcal{U}}$\fi}
\newcommand{\T}{\ifmmode{\mathcal{T}}\else${\mathcal{T}}$\fi}
\newcommand{\FF}{\ifmmode{\mathcal{F}}\else${\mathcal{F}}$\fi}
\newcommand{\LL}{\ifmmode{\mathcal{L}}\else${\mathcal{L}}$\fi}

\newtheorem{theorem}{Theorem}[section]
\newtheorem{proposition}[theorem]{Proposition}
\newtheorem{corollary}[theorem]{Corollary}
\newtheorem{lemma}[theorem]{Lemma}
\newtheorem{definition}[theorem]{Definition}
\newtheorem{remark}[theorem]{Remark}
\newtheorem{example}[theorem]{Example}

\title[Can one classify finite Postnikov pieces?]{Can one classify finite Postnikov pieces?}

\author{Jesper M. M{\o}ller}
\author{J\'er\^{o}me Scherer}

\thanks{The second author is supported by FEDER/MEC grant
MTM2007-61545 and by the MPI, Bonn.}


\begin{document}


\begin{abstract}
We compare the classical approach of constructing finite Postnikov
systems by $k$-invariants and the global approach of Dwyer, Kan, and
Smith. We concentrate on the case of $3$-stage Postnikov pieces and
provide examples where a classification is feasible. In general
though the computational difficulty of the global approach is
equivalent to that of the classical one.
\end{abstract}


\maketitle


\begin{quotation}
{\small \ldots all mathematics leads, doesn't it, sooner or later,
to some kind of human suffering.}
\end{quotation}
\begin{flushright}
   ``Against the Day", Thomas Pynchon
\end{flushright}

\section*{Introduction}
\label{sec intro}
Let $X$ be a finite Postnikov piece, i.e. a space with finitely many
non-trivial homotopy groups. Let us also assume for simplicity that
$X$ is simply connected. The classical theory of $k$-invariants
tells us that one can construct $X$ from Eilenberg-Mac\,Lane spaces
and a finite number of cohomology classes, the $k$-invariants, but
of course it might be difficult to compute them explicitly. This
computational difficulty is probably best illustrated by how
embarrassingly little one knows about the cohomology of Postnikov
pieces which are not $H$-spaces, see~\cite{MR0307242} for one of the
few examples where ``something" has been computed.

In \cite{dwks:towers}, Dwyer, Kan, and Smith propose a \emph{global}
approach. They provide in particular a model for the classifying
space of finite towers $X_n \rightarrow X_{n-1} \rightarrow \dots
\rightarrow X_1$ in which each fiber is a given Eilenberg-Mac\,Lane
space. We specialize to the case of $3$-stage Postnikov pieces, and
even further to fibrations of the form
$$
K(C, r) \times K(B, n) \rightarrow X \rightarrow K(A, m)
$$
with $1<m<n<r$. There exists a quite substantial literature about
this situation, let us mention especially Booth's work,
\cite{MR1817004} and \cite{MR2259270}, and Pave{\v{s}}i{\'c},
\cite{MR2175365}, \cite{MR2360776}. We explain first how the
Dwyer-Kan-Smith model provides a classifying space for such
fibrations and show in Corollary~\ref{cor:booth} that it coincides
with Booth's model from \cite{MR1817004}.

In the last section we compare then these two approaches, the
classical one based on $k$-invariants and the global one, and show
that they are basically equivalent. From the global point of view
what we must compute is a set of homotopy classes of lifts in a
fibration where the fiber is a product of Eilenberg-Mac\,Lane
spaces. It is quite remarkable how difficult it is to compute this,
compared to the elementary case when the fiber is a single
Eilenberg-Mac\,Lane space, a situation studied by the first author
in~\cite{MR910659}, and completely understood.

Consequently, a classification of three-stage Postnikov pieces will
be hopeless in general since it would necessitate the knowledge of
the cohomology of an arbitrary two-stage Postnikov piece. However,
classifications can be obtained in specific situations, and we
provide such examples along the way. The fact that we could not find
explicit computations in the literature motivated us to write this
note (one exception which we learned from Fernando Muro is Baues'
\cite[Section 10.5]{MR1404516}). Let us conclude with the comment
that this project grew up from a desire to understand the real scope
of the global classifying space model. Even though our conclusion
might seem rather pessimistic from a computational point of view, we
hope that the elegance of the global approach is still visible.
Moreover, our arguments do not exclude the possibility that
non-spectral sequence methods, applied to specific classes of
Postnikov pieces, may sometimes work better in the global context
than they do in the classical approach.

\medskip

\noindent {\bf Acknowledgements.} This work was achieved while the
second author was visiting the Max-Planck-Institut f\"ur Mathematik
in Bonn. We would like to thank Christian Ausoni and the referee for
helpful comments.

\section{Monoids of self-equivalences}
\label{sec:self-equivalences}

Let $X$ be a simply connected space. We consider some group-like
topological monoids consisting of (homotopy classes of) self
homotopy equivalences of $X$:
\begin{description}
\item[$\aut {}X$]  the topological monoid of self-homotopy
  equivalences of $X$,
\item[$\aut {*}X$]  the topological monoid of pointed self-homotopy
  equivalences of $X$,
\item[$\Aut {}X$]  the discrete group of components of $\aut{} X$.
\end{description}

In all cases the topological monoid structure is defined by
composition of maps. If $X$ happens to be an $H$-space, such as a
product of Eilenberg--Mac\,Lane spaces, then $\aut{}X$ also inherits
an $H$-space structure from $X$. These two structures are in general
not the same.

\begin{proposition}{\rm (\cite{MR0111028}, \cite{MR0154286},\cite[Chapter
IV]{MR1206474})} \label{prop:Baut}
There is a bijection of sets of homotopy classes of unpointed maps
$Y \rightarrow B \aut {}X$ and fiberwise homotopy types of
fibrations of the form $X \rightarrow E \rightarrow Y$.
\end{proposition}

If $t: Y \rightarrow B \aut {}X$ classifies such a fibration, one
often write $E = Y \times_t X$ for the total space and calls it a
twisted product. Much attention has been received by the set of
components $\Aut {}X$, but not so much by the space $\aut {}X$
itself. A nice exception is Farjoun and Zabrodsky's~\cite{MR546789}.

\section{Reminder on $2$-stage Postnikov systems}
\label{sec:2-stage}
In any introductory book on homotopy theory, such as
\cite[Chapter~IX]{MR516508}, one can read that a simply connected
space $E$ with only two non-trivial homotopy groups (say $\pi_m E
\cong A$ and $\pi_n E \cong B$ for $n > m$) is classified by a
$k$-invariant $k: K(A, m) \rightarrow K(B, n+1)$. This means that
$E$ has the homotopy type of the homotopy fiber of~$k$. How does
this relate to the approach described in the previous section?

We wish to understand the monoid $\aut {}{K(B, n)}$ and its
classifying space. From Proposition~\ref{prop:Baut} we infer that
two-stage Postnikov pieces $E$ with $\pi_m E \cong A$ and $\pi_n E
\cong B$ are in bijection with $[K(A, m), B\aut {}{K(B, n)}]$.

As a space $\aut {}{K(B, n)}$ is a product $\Aut {}B \times K(B,n)$;
this splitting is compatible with the $H$-space structure coming
from that of the Eilenberg-Mac\,Lane space, but not with the one we
are looking at, coming from composition. In fact $\aut {*}{K(B, n)}$
is weakly equivalent to the discrete monoid $\Aut {}{K(B, n)} \cong
\Aut {}B$. The weak equivalence is given by functoriality of the
$K(-, n)$ construction. Let us write $\varphi(\alpha)$ for the
pointed self-equivalence associated to the group
automorphism~$\alpha$. The map $\varphi$ splits the monoid map
$\pi_n: \aut {}{K(B, n)} \rightarrow \Aut {}B$. The fiber of $\pi_n$
over the identity is $\aut {1}{K(B, n)} \simeq {K(B, n)}$, on which
$\Aut {}B$ acts via $\varphi$ by conjugation. Thus we obtain a
description of the classifying space, see~\cite{MR0054250}.

\begin{lemma}
\label{lemma:pi*exact1}
The split exact sequence $K(B,n) \rightarrow \aut {}{K(B,n)}
\rightarrow \Aut {}B$ of topological monoids induces a split
fibration
\begin{equation*}
  \xymatrix@1{
   K(B,n+1)\ar[r] &
  {B\!\aut{}{K(B,n)}} \ar@<.5ex>[r] & {B\!\Aut {}B} \ar@<.5ex>[l]}
\end{equation*}
and thus ${B\!\aut{}{K(B,n)}}$ is the classifying space for
$(n+1)$-dimensional cohomology with local coefficients in~$B$.
\end{lemma}

\begin{proof}
The section given by functoriality of the construction of
Eilenberg-Mac\,Lane spaces is a map of monoids.
\end{proof}

We recover now the classification of fibrations with fibers $K(B,
n)$ as obtained by Steenrod in \cite[Section~23]{MR0298655} and Dold
\cite[Satz 12.15]{MR0198464}. They are classified by a single
$k$-invariant modulo the (non-trivial) action of $\Aut {}B$.

\begin{theorem}
\label{thm:2stageclassification}
The set of homotopy equivalent fibrations over a simply connected
space $Y$ with fibers $K(B, n)$ is in bijective correspondance with
$[Y, K(B, n+1)]/\Aut{}B \cong H^{m+1}(Y; B)/\Aut{}B$.
\end{theorem}

\begin{proof}
If we apply the functor $[Y,-]$ to the fibration from
Lemma~\ref{lemma:pi*exact1}, we obtain an exact sequence $[Y,
\Aut{}B]=\Aut{}B \rightarrow [Y, K(B,n+1)]=H^{n+1}(Y;B) \rightarrow
[Y, B\!\aut{}{K(B,n)}] \rightarrow *$ of sets and group actions.
\end{proof}

\begin{corollary}
\label{cor:2stageclassification}
Let $n > m > 1$ and $A, B$ be abelian groups. The set of $K(B,
n)$-fiber homotopy types over $K(A, m)$ is in bijection with
$H^{n+1}(K(A, m); B)/\Aut{}B$. \hfill{\qed}
\end{corollary}

Let us look at a basic example, which will serve as starting point
for examples of $3$-stage Postnikov pieces.

\begin{example}
\label{example E_2}
{\rm For $m=2$ and $n=3$, let us choose $A = B = \Z/2$ so
$\Aut{}{\Z/2} = 1$. Since $H^4(K(\Z/2, 2); \Z/2) \cong \Z/2$, there
are two homotopy spaces with the prescribed homotopy groups, namely
the product $K(\Z/2, 2) \times K(\Z/2, 3)$ and $E_2$ the homotopy
fiber of $Sq^2: K(\Z/2, 2) \rightarrow K(\Z/2, 4)$, the space
studied in \cite{MR0248811} by Milgram (and many others).
}
\end{example}

\section{The classical approach to $3$-stage Postnikov systems}
\label{sec naive}
In principle, the above theorem (and its corollary) can be used
inductively to classify $n$-stage Postnikov pieces. For a space $E$
with only three non-trivial homotopy groups $\pi_m E \cong A$,
$\pi_n E \cong B$, and $\pi_r E \cong C$ for $r> n > m >1$ we could
first construct $E[n]$, the $n$-th Postnikov section, which is given
by an element in $H^{m+1}(K(A, m); B)/\Aut{}B$ by
Corollary~\ref{cor:2stageclassification}. To reconstruct $E$ we will
then need to know the cohomology of $E[n]$, since the next
$k$-invariant lives in $H^{r+1}(E[n]; C)/\Aut{}C$. Our aim is to
study fiber homotopy types where the fiber is a product of two
Eilenberg-Mac\,Lane spaces.

Amazingly enough, we could not find a single example of
classification of $3$-stage Postnikov systems in the literature,
except \cite[Section~10.5]{MR1404516} where Baues deals with
$(n-1)$-connected $(n+2)$-types with $n \geq 4$. Let us treat
thoroughly one example, where we do the computations ``by hand". Its
interest also lies in the kind of computation one has to perform in
order to do the classification.

\begin{example}\label{exmp:mano}
{\rm Let us analyze fiber homotopy types of the form
\begin{equation*}
K(\Z/2, 5) \times K(\Z/2, 3) \to  E \to K(\Z/2, 2)
\end{equation*}
Thus $E$ has three non-trivial homotopy groups, all of them
isomorphic to $\Z/2$. There are two $k$-invariants. The first one
is a cohomology class $k_1 \in H^4(K(\Z/2, 2); \F_2) \cong \F_2
\langle Sq^2 \iota_2 \rangle$. Then the third Postnikov section
$E[3]$ is the homotopy fiber of $k_1$ and the second $k$-invariant
$k_2 \in H^6(E[3]; \F_2)$ is a class which restricts to zero in
$H^6(K(\Z/2, 3); \F_2)$ since we want the $3$-connected cover
$E\langle 3 \rangle$ to split as a product $K(\Z/2, 5) \times
K(\Z/2, 3)$.

When $k_1 = 0$, $k_2$ is a class in $H^6(K(\Z/2, 2) \times K(\Z/2,
3); \F_2)$ restricting to zero over $K(\Z/2, 3)$. By the K\"unneth
formula we see that $k_2$ lies in
$$
H^6(K(\Z/2, 2)) \oplus H^3(K(\Z/2, 2)) \otimes H^3(K(\Z/2, 3))
\oplus H^2(K(\Z/2, 2)) \otimes H^4(K(\Z/2, 3)).$$ There are thus 16
possible $k$-invariants, i.e. 16 different fiber homotopy types of
spaces $E$ over $K(\Z/2, 2) \times K(\Z/2, 3)$ with fiber $K(\Z/2,
5)$ such that $E\langle 3 \rangle \simeq K(\Z/2, 5) \times K(\Z/2,
3)$. This is not quite what we want. The group of components $\Z/2$
of $\aut {} {K(\Z/2, 2) \times K(\Z/2, 3)}$ acts on the 16
$k$-invariants by composition. It is easy to compute explicitly this
action of $\Z/2$: It acts trivially on 8 classes and identifies 4
pairs, so that we are left with 12 fiber homotopy types over
$K(\Z/2, 2)$.

When $k_1 = Sq^2 \iota_2$, let us denote by $E_2$ the homotopy
fiber. The mod $2$ cohomology of this space has been computed by
Milgram, \cite{MR0248811}, or Kristensen and Pedersen,
\cite{MR0324700}. It is an elementary Serre spectral sequence (for
the fibration $K(\Z/2, 3) \rightarrow E_2 \rightarrow K(\Z/2, 2)$)
argument to compute it in low degrees. We denote by $\iota_n$ the
non-zero class in $H^n(K(\Z/2, n); \F_2)$. In total degree $6$, the
only elements that survive are on the vertical axis -- $H^6(K(\Z/2,
3); \F_2)$ -- and the $\iota_2 \otimes Sq^1 \iota_3$ in bidegree
$(2, 4)$.

As the second $k$-invariant is a class in $H^6(E_2)$ restricting
to zero over $K(\Z/2, 3)$, it must be either zero or the class
corresponding to $\iota_2 \otimes Sq^1 \iota_3$. There are thus
only 2 fiber homotopy types over $E_2$. Now, in principle, there
could be an action of the group of self-equivalences of $E_2$
(isomorphic to $\Z/2$) on these two $k$-invariants, but as it
fixes zero, this action must be trivial. We have therefore also two
fiber homotopy types over $K(\Z/2, 2)$ covering $Sq^2 \iota_2$.}
\end{example}

%

The point of the example is that it illustrates well that one needs
to know the cohomology in low degrees of certain $2$-stage Postnikov
systems (and then identify the action of a group of
self-equivalences). This was easy here, but imagine the situation if
one would wish to compute fiber homotopy type over $K(\Z/2, 2)$ with
fiber $K(\Z/2, 3) \times K(\Z/2, 1000)$, or worse, to obtain a
classification in cases where the first $k$-invariant is not
primitive (say Schochet's \cite{MR0307242} homotopy fiber of the map
$K(\Z/2 \oplus \Z/2, 2) \rightarrow K(\Z/2, 4)$, represented in
cohomology by the product of the fundamental classes)!

\section{Spaces of lifts}
\label{sec:fibrations}
In this section we recall briefly the description and notation of
certain spaces of lifts from the work of the first author in
\cite{MR910659}. It deals with the case when the fiber is a single
Eilenberg-Mac\,Lane space. We then set up a spectral sequence to
treat the case of a Postnikov piece. Even in the case when the fiber
is a product of Eilenberg-Mac\,Lane spaces the description becomes
quickly complex. We start with some generalities about spaces of
lifts. Let us fix a fibration $p: Y \rightarrow Z$ and a map \func
uXY.

\begin{definition}\label{def:lifts}
{\rm The fiber containing $u \in \map {}XY$ of the induced fibration
$\underline{p}: \map {}XY \rightarrow \map {}XZ$ is the \emph{space
of lifts} $\Map uX{\emptyset}YZ = \{ v \in \map{}XY \mid pv=pu \}$
of all maps lying over $pu$.}
\end{definition}

Let \func{p_*}{[X,Y]}{[X,Z]} be the induced map of sets of homotopy
classes of maps.

\begin{lemma}\label{lemma:[X,Y]}
One has $[X,Y] \cong \displaystyle \coprod_{p_*u \in p_*[X,Y]} \pi_0
(\Map uX{\emptyset}YZ)/
  \pi_1(\map {}XZ,pu)$.
\end{lemma}
\begin{proof}
There are fibrations $\map{u}{X}{Y,Z} \rightarrow
\map{p_*^{-1}(p_*u)}XY \rightarrow \map {p_*u}XZ$
where $p_*u$ runs through the set $p_*[X,Y] \subset [X,Z]$.
\end{proof}

In the associated action
\begin{equation}\label{eq:action}
  \xymatrix@1{
  {\pi_1(\map {}XZ,pu) \times \pi_0(\Map uX{\emptyset}YZ)}
  \ar[r] &
  {\pi_0(\Map uX{\emptyset}YZ)}}
\end{equation}
the effect of an element $[h] \in \pi_1(\map{}XZ,pu)$ of the
fundamental group of the base space on the fibre $\Map
uX{\emptyset}YZ$ is given by \func{\overline{h}}{\Map
  uX{\emptyset}YZ}{\Map uX{\emptyset}YZ} where $\overline{h}$ is a
lift
\begin{equation*}
  \xymatrix{
    \{ 0 \} \times \Map uX{\emptyset}YZ \ar@{^(->}[rr] \ar@{^(->}[d] &&
    {\map {}XY} \ar[d]^-{\underline{p}} \\
    I \times  \Map uX{\emptyset}YZ \ar[r]^-{\mathrm{pr}_1}
                                   \ar@{.>}[urr]^-{\overline{h}} &
    I \ar[r]^-h & {\map{}XZ}}
\end{equation*}
of the homotopy $h \circ \mathrm{pr}_1$. Equivalently,
$\overline{h}$ is a solution
to the adjoint homotopy lifting problem $I \times  \Map
uX{\emptyset}YZ \times X \rightarrow Y$. Thus $\overline{h}$ is a
homotopy from the evaluation map $\overline{h}(0,v,x)=v(x)$ such
that $p\overline{h}(t,v,x)=h(t,x)$ is the given self-homotopy of
\func{pu}XZ. The end value of $\overline{h}$ takes $\Map
uX{\emptyset}YZ$ to itself.

\medskip

Assume now that the fibre of the fibration \func pYZ is the
Eilenberg--Mac\,Lane space $K(A,n)$. The primary difference between
the two lifts
\begin{equation*}
  \xymatrix@C=45pt{
    && Y \ar[d]^p \\
    {\Map  uX{\emptyset}YZ} \times X
    \ar@<.5ex>[urr]^-(.7){\mathrm{ev}}
    \ar@<-.5ex>[urr]_-(.7){u \circ \mathrm{pr}_2}
    \ar[r]^-{\mathrm{pr}_2} &
    X \ar[r]^{pu} & Z}
\end{equation*}
is an element of $\delta^n(\mathrm{ev},u \circ \mathrm{pr}_2)$ in
the group $H^n(\Map uX{\emptyset}YZ \times X;A)$. Let $\delta_i$ be
the components in $\prod H^i(\Map uX{\emptyset}YZ;H^{n-i}(X;A))$ of
$\delta^n(\mathrm{ev},u \circ \mathrm{pr}_2)$ under the isomorphism
\begin{equation*}
   H^n(\Map
  uX{\emptyset}YZ \times X;A)
  \cong \prod_{0 \leq i \leq n} H^i(\Map
  uX{\emptyset}YZ;H^{n-i}(X;A))
\end{equation*}
for the cohomology of a product. We can now state the
generalization of Thom's result, \cite{MR0089408}, obtained by the
first author.

\begin{theorem}\label{thm:thom}{\rm (M{\o}ller, \cite{MR910659})}
The map   \func{\prod \delta_i}{\Map uX{\emptyset}YZ} {\prod_{0 \leq
i \leq n} K(H^{n-i}(X;A),i)} is a homotopy equivalence.
\end{theorem}

In particular, $\pi_0(\Map uX{\emptyset}YZ) \cong H^n(X;A)$ and the
action \eqref{eq:action} takes the form of an action
\begin{equation*}
  \pi_1(\map{}XZ,pu) \times H^n(X;A) \to H^n(X;A)
\end{equation*}
of the group $\pi_1(\map{}XZ,pu)$ on the {\em set\/} $H^n(X;A)$. How
can we describe this action?

\begin{lemma}
Let $\mathrm{ev} \in H^n(\map {}X{K(A,n)} \times X;A)$ be the
evaluation map. Write $\mathrm{ev}= \sum \mathrm{ev}_i$ as a sum of
cohomology classes under the K\"unneth isomorphism
  \begin{equation*}
    H^n(\map {}X{K(A,n)} \times X;A ) \cong
   \bigoplus_{i+j=n} H^i(\map {}X{K(A,n)};H^j(X;A))
  \end{equation*}
  Then
  \begin{equation*}
     \prod \mathrm{ev}_i \colon
     \map {}X{K(A,n)} \to \prod_{i+j=n} K(H^j(X;A),i)
  \end{equation*}
is a homotopy equivalence. \hfill{\qed}
\end{lemma}

We are now ready for the promised spectral sequence computing the
homotopy groups of the space of lifts in a fibration where the fiber
has more than a single non-trivial homotopy group (there is an
analogous spectral sequence when the source $X$ is a finite
CW-complex). It is obtained by decomposing the fiber by its
Postnikov sections. The case of a space of sections has been studied
in great detail by Legrand, \cite{MR668016}. It generalizes work of
Shih, \cite{MR0144348}, on limited and non-abelian spectral
sequences. The bigrading we have chosen here agrees with that in
\cite[Theorem~5.3]{MR1038731}.

\begin{corollary}
\label{cor:liftSS}
Suppose that $F \to Y \to Z$ is a split fibration where the fibre
$F$ is a finite Postnikov piece, connected and simple. Let $u: X
\rightarrow Z$ be a map. Then there is third octant homology
spectral sequence ($i+j \geq 0$ and $i \leq 0$)
  \begin{equation*}
    E^2_{ij} = H^{-i}(X;\pi_{j}(F)) \Longrightarrow \pi_{i+j}(\Map
    uX{\emptyset}YZ)
  \end{equation*}
converging to the homotopy groups of the space of lifts.
\hfill{\qed}
\end{corollary}

In principle, the cohomology groups appearing in the spectral
sequence are to be understood with local coefficients defined by the
choice of a lift. The space of lifts here is not empty since we
assume for simplicity that the fibration has a section. The case
when the fiber has two non-trivial homotopy groups is already
interesting.

\begin{example}\label{exmp:Fisaproduct}
{\rm Suppose that the fibre $F=K(A,m) \times K(B,n)$ with $m < n$.
In that case the spectral sequence is concentrated on two lines and
yields a long exact sequence. It can be identified with the homotopy
long exact sequence of the fibration
  \begin{equation*}
    \Map uZ{\emptyset}{Y}{\overline Y[m]} \longrightarrow
    \Map uZ{\emptyset}{Y}{Z} \longrightarrow
    \Map uZ{\emptyset}{\overline Y[m]}{Z}
  \end{equation*}
where $\overline Y[m]$ denotes the fiberwise Postnikov section, i.e.
the map $Y \rightarrow Z$ factors through $\overline Y[m]$ and the
homotopy fiber of $\overline Y[m] \rightarrow Z$ is $F[m] = K(A,m)$.
We deduce from Theorem~\ref{thm:thom} that $\Map
uZ{\emptyset}{Y}{\overline Y[m]} \simeq \prod K(H^{n-i}(Z;B),i)$ and
$\Map uZ{\emptyset}{\overline Y[m]}{Z} \simeq \prod
K(H^{m-i}(Z;A),i)$. Hence the long exact sequence terminates in
particular with
  \begin{equation*}
    H^{m-1}(X;A) \to  H^{n}(X;B) \to \pi_0 \Map
    uZ{\emptyset}{Y}{Z} \to  H^{m}(X;A)
  \end{equation*}
Note that even though the fibre is a product, the $k$-invariant
$Y[m]  \to K(B,n+1)$ may not be trivial (it only restricts to $0$ on
the fibre) and therefore the $k$-invariant of the above fibration
may not be trivial either so that the sequence does not split!}
\end{example}

This indicates that, as soon as there are more than one
non-trivial homotopy group in the fiber, it will be difficult even
to compute the number of homotopy classes of lifts, in contrast
with Theorem~\ref{thm:thom}.

\section{The Dwyer-Kan-Smith model}
\label{sec DKS}
Let us now look at the ``global" point of view on Postnikov pieces.
Instead of adding iteratively one Eilenberg-Mac\,Lane space at a
time, one can also try to understand how to add \emph{all} homotopy
groups at once. This is the approach followed by Dwyer, Kan, and
Smith in \cite{dwks:towers}. In this section we will see how it
specializes to the case of $3$-stage Postnikov pieces and which
modifications we need to obtain explicit classification results.

Let $G$ be a space and consider the functor $\Phi$ which sends an
object of $Spaces \downarrow B \aut{} {G}$, i.e. a map $t: X
\rightarrow B \aut{} {G}$, to the twisted product $X \times_t G$,
see Section~\ref{sec:self-equivalences}. Dwyer, Kan, and Smith
describe a right adjoint $\Psi$ in \cite[Section~4]{dwks:towers}.
They find first a model for $\aut{} {G}$ which is a (simplicial)
group and thus acts on the left on $\map{} {G} Z$ for any space $Z$.
This induces a map $r: B \aut{} {G} \rightarrow B\aut{}{\map{} {G}
Z}$. The functor $\Psi$ then sends $Z$ to the projection map from
the twisted product, $B \aut{} {G} \times_r \map{} {G} Z \rightarrow
B \aut{} {G}$. This allows us to immediately construct a classifying
space for towers, in our case they will be of length~$2$.

\begin{theorem}\label{thm:dks1}{\rm (Dwyer, Kan, Smith,
\cite{dwks:towers})}
The classifying space for towers of the form $Z \xrightarrow{q} Y
\xrightarrow{p} X$, where the homotopy fiber of $p$ is $G$ and that
of $q$ is $H$, is $B \aut{} {G} \times_r \map{} {G} {B \aut{} {H}}$.
\end{theorem}

Fix now a fibration $H \rightarrow F \rightarrow G$ where we think
about the spaces $H$ and $G$ as simpler, in particular the spaces
$\aut{}H$ and $\aut{}G$ should be accessible. Such a fibration is
classified by a map $s: G \rightarrow B \aut{} H$ and so $F$ is the
\emph{twisted product} $G \times_s H$. To construct $B \aut{} {F}$,
one simply needs to refine a little the analysis done by Dwyer, Kan,
and Smith. Let us denote by $\map{[s]} {G} {B \aut{} {H}}$ the
components of the mapping space corresponding to the orbit of the
map $s$ defined above under the action of $\Aut{} {G}$.

\begin{lemma}\label{lem:dks}
Let $H \rightarrow F \rightarrow G$ be any fibration, classified by
a map $s: G \rightarrow B \aut{} H$.  The space $B \aut{} {G}
\times_r \map{[s]} {G} {B \aut{} {H}}$ classifies towers $Z
\xrightarrow{\beta} Y \xrightarrow{\alpha} X$ where the homotopy
fiber of $\alpha$ is $G$, that of $\beta$ is $H$, and that of the
composite $\alpha \circ \beta$ is~$F$.
\end{lemma}

\begin{proof}
Since $B \aut{} {G} \times_t \map{[s]} {G} {B \aut{} {H}}$ is a
subspace of the classifying space for towers $Z \rightarrow Y
\rightarrow X$ over $X$ with fibers $G$ and $H$, it classifies some
of them. We claim that the fiber of the composite map $Z \rightarrow
X$ is precisely $F$.

From the adjunction property a map $X \rightarrow B \aut{} {G}
\times_r \map{[s]} {G} {B \aut{} {H}}$ corresponds to a map $t': X
\times_t G \rightarrow B \aut{} {H}$, which yields a space $E = X
\times_t G \times_{t'} H$. The fiber we must identify is thus the
homotopy pull-back of the diagram $E \rightarrow X \times_t G
\leftarrow G$. In other words it is the twisted product
corresponding to the composite map $G \rightarrow X \times_t G
\rightarrow B\aut{}H$, which is homotopic to~$s$. This means that
the homotopy fiber is~$F$.
\end{proof}

To find an description of $B\aut{} F$ in terms of $G$ and $H$ is a
more difficult task, because in general not all fibrations with
fiber $F$ come from a tower as above. However there are situations
where this is so. Let us consider a \emph{homotopy localization}
functor $L$, like Postnikov sections, Quillen plus-construction, or
localization at a set of primes, see \cite{Dror}. What matters for
us is that there are natural maps $\eta: X \rightarrow LX$ for all
spaces $X$, and that $L$ sends weak equivalences to weak
equivalences.

\begin{theorem}\label{thm:dks}
Let $L$ be a homotopy localization functor and consider a fibration
$\bar L F \rightarrow F \xrightarrow{\eta} LF$, classified by a map
$s: LF \rightarrow B \aut{} {\bar L F}$. Then $B\aut{}F$ is $B
\aut{} {LF} \times_r \map{[s]} {LF} {B \aut{} {\bar L F}}$.
\end{theorem}

\begin{proof}
Let $F \rightarrow Z \rightarrow X$ be any fibration over $X$. It is
possible to construct a fiberwise version of $L$, i.e. obtain a new
fibration $LF \rightarrow Y \rightarrow X$ such that the diagram
\[
\diagram
  F \ar[r] \ar[d]_\eta & Z \ar[r] \ar[d] & X \ar@{=}[d]\\
  LF \ar[r] & Y \ar[r] & X
\enddiagram
\]
commutes, \cite[Theorem~F.3]{Dror}. So any fibration comes from a
tower $Z \rightarrow Y \rightarrow X$. In particular this
construction can be applied to the universal fibration $F
\rightarrow B\aut*F \rightarrow B\aut{}F$ and this yields a map
$B\aut{}F \rightarrow B \aut{} {LF} \times_r \map{} {LF} {B \aut{}
{\bar L F}}$, which factors through the component $B \aut{} {LF}
\times_r \map{[s]} {LF} {B \aut{} {\bar LF}}$ by
Lemma~\ref{lem:dks}. There is a forgetful map going the other way,
and both composites are homotopic to the identity by uniqueness of
the classifying space.
\end{proof}

We are mainly interested in $3$-stage Postnikov systems in this
note. Consider thus a $3$-stage Postnikov piece $E$ as being the
total space of a fibration of the form $F \rightarrow E \rightarrow
K(A, m)$. The fiber $F$ is a space with only two non-trivial
homotopy groups and the fibration is classified by a map $K(A, m)
\rightarrow B \aut{}F$, see Theorem~\ref{prop:Baut}. We understand
now the monoid of self-equivalences of a space with two non-trivial
homotopy groups.

\begin{corollary}\label{cor:booth}
Let $F$ be a simply connected $2$-stage Postnikov piece, with $\pi_m
F \cong A$, $\pi_n F \cong B$, $k$-invariant $k$, and $n>m$. Then
$B\aut{}F$ is $B \aut{} {K(A, m)} \times_r \map{[k]} {K(A, m)} {B
\aut{} {K(B, n)}}$.
\end{corollary}

\begin{proof}
The $m$-th Postnikov section $F \rightarrow F[m]$ is a homotopy
localization functor.
\end{proof}

Let us specialize even further, and assume that the $k$-invariant is
trivial, that is, we are looking at a fiber which is a product of
two Eilenberg-Mac\,Lane spaces. Such a model has been independently
constructed by Booth in \cite{MR1817004}.

\begin{corollary}\label{cor:booth2}
Let $A$ and $B$ be two abelian groups and $n > m$. Then $B \aut{}
{K(A, m) \times K(B, n)} \simeq B \aut{} {K(A, m)} \times_r \map{c}
{K(A, m)} {B \aut{} {K(B, n)}}$, where $c$ is the constant map. The
projection $B \aut{} {K(A, m) \times K(B, n)} \rightarrow B \aut{}
{K(A, m)}$ has a section.
\end{corollary}

\begin{proof}
The orbit of the constant map is reduced to the constant map.
\end{proof}

The computation of the set of components of $\aut {} {K(A, m) \times
K(B, n)}$ is straightforward, compare with Shih's \cite{MR0160209},
or the matrix presentation used in \cite[Section~1]{MR2259270}.

\begin{corollary}\label{cor:setofcomponents}
Let $A$ and $B$ be two abelian groups and $n > m>1$ be integers.
Then the group $\Aut{} {K(A, m) \times K(B, n)}$ is a split
extension of $\Aut{} A \times \Aut{} B$ by $H^n(K(A, m); B)$.
\hfill{\qed}
\end{corollary}

\section{Comparing the classical with the global approach}
\label{sec compare}
The classical approach to finite $n$-stage Postnikov pieces goes
through the computation of the cohomology of a $(n-1)$-stage
Postnikov piece. This is theoretically feasible via a Serre spectral
sequence computation, but practically very hard because of the
differentials. What about the global approach?

We consider the case of fiber homotopy types over $K(A, m)$ with
fiber $K(B, n) \times K(C, r)$ with $1<m<n<r$ as before. In
principle we only need to compute the set of homotopy classes $[K(A,
m), B \aut{} {(K(B, n) \times K(C, r)}]$ and we have a model for
this classifying space. The only sensible way we could think of to
compute this is by using the split fibration
$$
\map{c}{K(B,
n)}{B\aut{}{K(C, r)}} \rightarrow B \aut{} {(K(B, n) \times K(C, r)}
\rightarrow B\aut{}{K(B, n)}
$$
obtained in Corollary~\ref{cor:booth2}. Thus for each first
$k$-invariant $k_1: K(A, m) \rightarrow K(B, n+1)$ we must
understand the set of components of the space of lifts into $B\aut{}
{(K(B, n) \times K(C, r)}$.

\begin{example}\label{exmp:dks}
{\rm Let us again analyze fiber homotopy types of the form
\begin{equation*}
K(\Z/2, 5) \times K(\Z/2, 3) \to  E \to K(\Z/2, 2)
\end{equation*}
We will now do the computation \emph{globally}. Let us write shortly
$K_n$ for the space $K(\Z/2, n)$. The classifying space is $K_4
\times_t \map{c}{K_3} {K_6}$. Consider now the sectioned fibration
\begin{equation*}
  \xymatrix@1{
  K_6 \times K_3 \times K_2 \times K_1 =
  \map c{K_3}{K_6}
  \ar[r] &
 {B\!\aut{}{K_3 \times K_5}} \ar@<.5ex>[r] &
 {B\!\aut{}{K_3} = K_4} \ar@<.5ex>[l]^s}
\end{equation*}
so that $[K_2,B\!\aut{}{K_3 \times K_5}]$ is the disjoint union of
the components of $\map {}{K_2}{B\!\aut{}{K_3 \times K_5}}$ which
lie over $0$ and those which lie over $Sq^2 \iota_2$ in
$\map{}{K_2}{K_4}$. By Lemma~\ref{lemma:[X,Y]} these two sets can be
computed as quotients of sets of components of spaces of lifts under
the action of a fundamental group.

Let us do that. Over zero, there is no mystery, the space of lifts
is $\map {}{K_2}{K_6 \times K_3 \times K_2 \times K_1}$ and the
fundamental group in question is $\pi_1 \map {}{K_2}{K_4} \cong
\Z/2$. It is straightforward to see that the 16 components of the
mapping space are grouped in 12 orbits. Over $Sq^2\iota_2$, we are
looking at the space of lifts as in the following diagram:

\[\xymatrix@1{
& K_4 \times_t \map{c}{K_3} {K_6} \ar[d] \cr K_2 \ar[r]^{Sq^2}
\ar[ur] & K_4}
\]
This is equivalent by the Dwyer-Kan-Smith adjunction
\cite[Section~4]{dwks:towers} to the subspace of maps
$\map{}{E_2}{K_6}$ which restrict trivially to~$K_3$. From the 16
possible components we are left with~two, compare with
Example~\ref{exmp:Fisaproduct}. The action of $\pi_1 \map
{Sq^2}{K_2}{K_4} \cong \Z/2$ is trivial and it seems we have redone
here as well the same computation as in Example~\ref{exmp:mano}.}
\end{example}

Let us carefully check whether we have really redone the same
computations as in the classical approach.

Our typical study case is that of a space with three non-trivial
homotopy groups $A$, $B$, and $C$, in degree respectively $m$, $n$,
and $r$, with $1<m<n<r$. In the classical approach we use for each
possible first $k$-invariant $k_1: K(A, m) \rightarrow K(B, n+1)$
the corresponding Serre spectral sequence $H^p(K(A, m); H^q(K(B, n);
C))$ of which we only need the $p+q=r+1$-diagonal to determine the
possible values of the second $k$-invariant.

In the global approach we wish to compute, for each possible first
$k$-invariant $k_1: K(A, m) \rightarrow K(B, n+1)$, the set of
components of the space of lifts indicated in the diagram

\[\xymatrix@1{
& K(B, n+1) \times_t \map{c}{K(B, n)}{K(C, r+1)} \ar[d] \cr K(A, m)
\ar[r]^{k_1} \ar[ur] & K(B, n+1)}
\]
Since the mapping space $\map{c}{K(B, n)}{K(C, r+1)}$ is a product
of Eilenberg-Mac\,Lane spaces $K(H^{n+1-j}(K(B, n); C), j)$ for $1
\leq j \leq n+1$, Corollary~\ref{cor:liftSS} yields a spectral
sequence of the form $E_2^{ij}=H^{-i}(K(A, m); H^{r+1-j}(K(B, n);
C))$ with differential $d_2$ of bidegree $(-2, 1)$.

Techniques due to Legrand, \cite{MR668016} and \cite{MR0500962},
allow to prove the following result.

\begin{proposition}\label{prop:compareSS}{\rm (Didierjean and Legrand,
\cite[Th\'eor\`eme~2.2]{MR746503})}
Suppose that $F \to Y \to Z$ is a fibration where the fibre $F$ is a
connected, finite Postnikov piece. The spectral sequences converging
to the homotopy groups of the space of lifts $\Map uX{\emptyset}YZ$
defined from the skeletal filtration of $X$, and the one defined by
the Postnikov decomposition of $F$ are isomorphic.
\end{proposition}

\begin{proof}
The same argument as in \cite{MR746503} for spaces of sections
applies for spaces of lifts. It relies on the techniques developed
in~\cite{MR668016}. Alternatively one could identify the space of
lifts as a space of sections (of the pull-backed fibration) and
apply directly Didierjean and Legrand's result.
\end{proof}

\begin{remark}
\label{rem:historic}
{\rm This kind of spectral sequence appeared maybe first in work of
Federer, \cite{MR0079265}. It also appears in Switzer,
\cite{MR638816}, in both forms, but he does not compare them
however. When the target $Y$ is a spectrum rather than a space, the
spectral sequences are the Atiyah-Hirzebruch one and the Postnikov
one. Maunder proved they coincide, \cite{MR0150765}. When $Y$ is a
space, like here, cosimplicial technology allowed Bousfield to
construct such spectral sequences yet in another way,
\cite{MR1017155}.}
\end{remark}

We now come back to our Postnikov pieces. The above proposition
allows us to identify the spectral sequence coming from a Postnikov
decomposition of $\map{c}{K(B, n)}{K(C, r+1)}$ with the one coming
from the skeletal filtration of $K(A, m)$. The last step is to
identify this second spectral sequence. Let us first regrade the
spectral sequence by setting $p = -i$ and $q= r+1-j$, so our
$E_2$-term looks like $E_2^{pq} = H^p(K(A, m); H^{q}(K(B, n); C))$
(and the differential $d_2$ has bidegree $(2, -1)$). This spectral
sequence is concentrated in the first quadrant, in the horizontal
band $0 \leq q \leq r+1$. It converges to $\pi_{p+q-r-1} \Map
{k_1}{K(A, m)}{\emptyset}{K(B, n+1)}{K(B, n+1) \times_t \map{c}{K(B,
n)}{K(C, r+1)}}$.

\begin{theorem}\label{prop:compareSS2}
Let $r>n>m>1$ be integers and $A, B, C$ be abelian groups. For any
$k$-invariant $k_1: K(A, m) \rightarrow K(B, n+1)$, the part of the
Postnikov spectral sequence concentrated in degrees $p+q \leq r+1$
computing the homotopy groups of the space of lifts into $K(B, n+1)
\times_t \map{c}{K(B, n)}{K(C, r+1)}$ over $k_1$ is isomorphic to
the corresponding part of the cohomological Serre spectral sequence
(with coefficients in $C$) for the fibration $K(B, n) \rightarrow
K(A, m) \times_{k_1} K(B, n) \rightarrow K(A, m)$.
\end{theorem}

\begin{proof}
The fiber $\map{c}{K(B, n)}{K(C, r+1)}$ is a connected and finite
Postnikov piece with abelian fundamental group, so that the
Postnikov spectral sequence from Corollary~\ref{cor:liftSS} exists.
From the previous proposition we know that it actually coincides
with the spectral sequence defined by the skeletal filtration of
$K(A, m)$.

Instead of looking at the $E_2$-term we will work with the
$E_1$-term. We write $K(A, m)_k \subset K(A, m)$ for the $k$-th
skeleton, and $p: K(A, m) \times_{k_1} K(B, n) \rightarrow K(A, m)$
for the natural projection. By the Dwyer-Kan-Smith adjunction lifts
over $K(A, m)_k$  correspond to maps from the preimage under $p$ to
$K(C, r+1)$. This is precisely the filtration in the Serre spectral
sequence. All differentials in the triangle $p+q \leq r+1$ remain in
the band $0 \leq q \leq r+1$, in which the $E_2$-term of the
Postnikov sequence is abstractly isomorphic to the $E_2$-term of the
Serre spectral sequence thanks to the regrading we have performed
(for $q > r+1$ it is zero).
\end{proof}

\begin{remark}\label{rem:compareSS}
{\rm Let $r>n>m>1$ be integers and $A, B, C$ be abelian groups. We
have seen two approaches to compute the number of fiber homotopy
types $X$ over $K(A, m)$ with fiber $K(B, n) \times K(C, r)$ such
that $X[n]$ is classified by a given $k$-invariant $k_1: K(A, m)
\rightarrow K(B, n+1)$. The one we have called the global one
computes the set of components of a space of lifts via a Postnikov
spectral sequence. Since the diagonal $p+q-r-1=0$ is contained (as
the edge) in the triangle we have been able to analyze in
Theorem~\ref{prop:compareSS2}, we see that this computation is
\emph{exactly the same} as the classical one, where one is looking
for the second $k$-invariant.}
\end{remark}

%
%


\bibliographystyle{amsplain}

\begin{thebibliography}{10}

\bibitem{MR0111028}
M.~G. Barratt, V.~K. A.~M. Gugenheim, and J.~C. Moore, \emph{On
semisimplicial
  fibre-bundles}, Amer. J. Math. \textbf{81} (1959), 639--657.

\bibitem{MR1404516}
H.-J. Baues, \emph{Homotopy type and homology}, Oxford Mathematical
Monographs,
  The Clarendon Press Oxford University Press, New York, 1996, Oxford Science
  Publications.

\bibitem{MR1817004}
P.~I. Booth, \emph{Fibrations with product of
{E}ilenberg-{M}ac{L}ane space
  fibres. {I}}, Groups of homotopy self-equivalences and related topics
  (Gargnano, 1999), Contemp. Math., vol. 274, Amer. Math. Soc., Providence, RI,
  2001, pp.~79--104.

\bibitem{MR2259270}
\bysame, \emph{An explicit classification of three-stage {P}ostnikov
towers},
  Homology, Homotopy Appl. \textbf{8} (2006), no.~2, 133--155 (electronic).

\bibitem{MR1017155}
A.~K. Bousfield, \emph{Homotopy spectral sequences and
obstructions}, Israel J.
  Math. \textbf{66} (1989), no.~1-3, 54--104.

\bibitem{MR746503}
A.~Didierjean and A.~Legrand, \emph{Suites spectrales de {S}erre en
homotopie},
  Ann. Inst. Fourier (Grenoble) \textbf{34} (1984), no.~2, 227--242.

\bibitem{MR0198464}
A.~Dold, \emph{Halbexakte {H}omotopiefunktoren}, Lecture Notes in
Mathematics,
  vol.~12, Springer-Verlag, Berlin, 1966.

\bibitem{MR546789}
E.~Dror and A.~Zabrodsky, \emph{Unipotency and nilpotency in
homotopy
  equivalences}, Topology \textbf{18} (1979), no.~3, 187--197.

\bibitem{dwks:towers}
W.~G. Dwyer, D.~M. Kan, and J.~H. Smith, \emph{Towers of fibrations
and
  homotopical wreath products}, J. Pure Appl. Algebra \textbf{56} (1989),
  no.~1, 9--28.

\bibitem{Dror}
E.~Dror Farjoun, \emph{Cellular spaces, null spaces and homotopy
localization},
  Lecture Notes in Mathematics, vol. 1622, Springer-Verlag, Berlin, 1996.

\bibitem{MR0079265}
H.~Federer, \emph{A study of function spaces by spectral sequences},
Trans.
  Amer. Math. Soc. \textbf{82} (1956), 340--361.

\bibitem{MR0324700}
L.~Kristensen and E.~K. Pedersen, \emph{The {${\mathcal A}$}-module
structure
  for the cohomology of two-stage spaces}, Math. Scand. \textbf{30} (1972),
  95--106.

\bibitem{MR0500962}
A.~Legrand, \emph{Sur les groupes d'homotopie des sections continues
d'une
  fibration en groupes}, C. R. Acad. Sci. Paris S\'er. A-B \textbf{286} (1978),
  no.~20, A881--A883.

\bibitem{MR668016}
\bysame, \emph{Homotopie des espaces de sections}, Lecture Notes in
  Mathematics, vol. 941, Springer-Verlag, Berlin, 1982.

\bibitem{MR0150765}
C.~R.~F. Maunder, \emph{The spectral sequence of an extraordinary
cohomology
  theory}, Proc. Cambridge Philos. Soc. \textbf{59} (1963), 567--574.

\bibitem{MR1206474}
J.~P. May, \emph{Simplicial objects in algebraic topology}, Chicago
Lectures in
  Mathematics, University of Chicago Press, Chicago, IL, 1992, Reprint of the
  1967 original.

\bibitem{MR0248811}
R.~James Milgram, \emph{The structure over the {S}teenrod algebra of
some
  {$2$}-stage {P}ostnikov systems}, Quart. J. Math. Oxford Ser. (2) \textbf{20}
  (1969), 161--169.

\bibitem{MR910659}
J.~M. M{\o}ller, \emph{Spaces of sections of {E}ilenberg-{M}ac
{L}ane
  fibrations}, Pacific J. Math. \textbf{130} (1987), no.~1, 171--186.

\bibitem{MR1038731}
\bysame, \emph{On equivariant function spaces}, Pacific J. Math.
\textbf{142}
  (1990), no.~1, 103--119.

\bibitem{MR0054250}
P.~Olum, \emph{On mappings into spaces in which certain homotopy
groups
  vanish}, Ann. of Math. (2) \textbf{57} (1953), 561--574.

\bibitem{MR2175365}
P.~Pave{\v{s}}i{\'c}, \emph{On the group {${\rm Aut}\sb {\#}(X\sb
  1\times\cdots\times X\sb n)$}}, Topology Appl. \textbf{153} (2005), no.~2-3,
  485--492.

\bibitem{MR2360776}
\bysame, \emph{Reducibility of self-homotopy equivalences}, Proc.
Roy. Soc.
  Edinburgh Sect. A \textbf{137} (2007), no.~2, 389--413.

\bibitem{MR0307242}
C.~Schochet, \emph{A two-stage {P}ostnikov system where
  {$E\sb{2}\not=E\sb{\infty }$} in the {E}ilenberg-{M}oore spectral sequence},
  Trans. Amer. Math. Soc. \textbf{157} (1971), 113--118.

\bibitem{MR0144348}
W.~Shih, \emph{Homologie des espaces fibr\'es}, Inst. Hautes
\'Etudes Sci.
  Publ. Math. (1962), no.~13, 88.

\bibitem{MR0160209}
\bysame, \emph{On the group {${\mathcal E}[X]$} of homotopy
equivalence maps},
  Bull. Amer. Math. Soc. \textbf{70} (1964), 361--365.

\bibitem{MR0154286}
J.~Stasheff, \emph{A classification theorem for fibre spaces},
Topology
  \textbf{2} (1963), 239--246.

\bibitem{MR0298655}
N.~E. Steenrod, \emph{Cohomology operations, and obstructions to
extending
  continuous functions}, Advances in Math. \textbf{8} (1972), 371--416.

\bibitem{MR638816}
R.~M. Switzer, \emph{Counting elements in homotopy sets}, Math. Z.
\textbf{178}
  (1981), no.~4, 527--554.

\bibitem{MR0089408}
R.~Thom, \emph{L'homologie des espaces fonctionnels}, Colloque de
topologie
  alg\'ebrique, Louvain, 1956, Georges Thone, Li\`ege, 1957, pp.~29--39.

\bibitem{MR516508}
G.~W. Whitehead, \emph{Elements of homotopy theory}, Graduate Texts
in
  Mathematics, vol.~61, Springer-Verlag, New York, 1978.

\end{thebibliography}
\providecommand{\bysame}{\leavevmode\hbox
to3em{\hrulefill}\thinspace}
\providecommand{\MR}{\relax\ifhmode\unskip\space\fi MR }
\providecommand{\MRhref}[2]{%
  \href{http://www.ams.org/mathscinet-getitem?mr=#1}{#2}
} \providecommand{\href}[2]{#2}



\bigskip
{\small
\begin{center}
\begin{minipage}[t]{7 cm}
Jer\^{o}me Scherer\\ Departament de Matem\`atiques,\\
Universitat Aut\`onoma de Barcelona,\\ E-08193 Bellaterra, Spain
\\\textit{E-mail:}\texttt{\,\,jscherer@mat.uab.es}
\end{minipage}
\begin{minipage}[t]{7 cm}
Jesper M. M{\o}ller \\  Matematisk Institut \\ Universitetsparken 5
\\DK-2100 K{\o}benhavn
\\\textit{E-mail:}\texttt{\,moller@math.ku.dk},
\end{minipage}
\end{center}}

\end{document}